\documentclass[12pt]{article}
\usepackage{amssymb}
\newtheorem{theorem}{Theorem}

\newtheorem{proposition}{Proposition}
\title{A Combinatorial Method for Counting Smooth Numbers in Sets of Integers}
\author{Ernie Croot}

\begin{document}
\maketitle
\begin{abstract}  
In this paper we prove a result for determining the number of integers without 
large prime factors lying in a given set $S$.  We will apply it to give an easy proof
that certain sufficiently dense sets $A$ and $B$ always produce the expected 
number of ``smooth'' sums $a+b$, $a \in A$, $b \in B$.  The proof of this result is 
completely combinatorial and elementary.
\end{abstract}

\section{Introduction}

Given a set $S$, a common question one tries to answer is whether
$S$ contains the expected number of ``$y$-smooth'' integers, which
are those integers having no prime divisors greater than $y$.  We denote
the number of integers in $S$ with this property by $\Psi(S,y)$; and for
a number $x > 0$, one denotes the set of all $y$-smooths positive integers
$\leq x$ by $\Psi(x,y)$.  So, 
$$
\Psi(\{1,2,..., \lfloor x \rfloor \}, y)\ =\ \Psi(x,y).
$$
If $S \subseteq \{1,2,...,x\}$, then, all things being equal, one would 
expect that
\begin{equation} \label{hope}
{\Psi(S,y) \over |S|}\ \sim\ {\Psi(x,y) \over x}.
\end{equation}
For example, fix a real number $0 < \theta \leq 1$ and an integer
$a \neq 0$, and let $S$ be the
set of numbers of the form $p+a$, where $p \leq x$ runs through the
primes; $S$ is often called a set of ``shifted primes''.  It is conjectured
that
\begin{equation} \label{difficult}
\Psi(S,x^\theta)\ \sim\ {\pi(x) \Psi(x,x^\theta) \over x}\ \sim\ \rho(\theta^{-1}) \pi(x),
\end{equation}
where 
$$
\rho(u)\ =\ \lim_{x\to \infty} {\Psi(x,x^{1/u}) \over x}.
$$
This function $\rho$ is called Dickman's function, and it was proved in
\cite{dickman} that the limit exists.   
Unfortunately, proving (\ref{difficult}) remains a difficult, open problem; however,
in \cite{friedlander}, J. B. Friedlander gave a beautiful proof that 
$\Psi(S,x^\theta) \gg \pi(x)$ for $\theta > (2 \sqrt{e})^{-1}$, and 
in \cite{baker}, R. Baker and G. Harman
proved that for $\theta \geq 0.2961$,
$$
\Psi(S,x^\theta)\ >\ {x \over \log^\alpha x},
$$
for some $\alpha > 1$ and $x > x_0(a)$.

There are several methods for attacking the general question of
proving that (\ref{hope}) holds for a particular set $S$, one such method
involves exponential sums and the circle method, and another uses a
Buchstab identity, in combination with a sieve method (such as the Large
Sieve).  

In this paper we offer a novel way of showing that sets $S$ have the expected
number of $x^\theta$-smooths, and the conditions that $S$ needs to satisfy,
in order for this method to work, are simpler than those required by 
other methods (such as an application of Buchstab's identity).
Before we can state what these conditions are, we first introduce the notion
of a ``Local-Global Set'', which we abbreviate as LG set:
\bigskip

\noindent {\bf Definition.  }  We say that  
$N \subseteq \{2,3,...,x\}$ is an LG set with parameters
$\epsilon$, $c$ and $x$ if and only if following two conditions hold:

1.  For any pair of distinct members $n_1,n_2 \in N$ we have 
lcm$(n_1,n_2) > x$;

2.  All but at most $\epsilon x$ of the integers $m \leq x$ are divisible by
some $n \leq x^c$ with $n \in N$.
\bigskip

\noindent {\bf Notes:}  From condition 1 we know that if $m \leq x$ is divisible by
one of these $n$'s, then this $n$ must be unique, else if 
$n_1,n_2 \in N$ are distinct and if $n_1 | m$ and $n_2 | m$, then
lcm$(n_1,n_2) | n$, which implies $n \geq {\rm lcm}(n_1,n_2) > x$, 
contradiction.  We also note that condition 2 implies
\begin{eqnarray}\label{N_sum}
\sum_{n \in N \atop n < x^c} 
{1 \over n}\ &=&\ {1 \over x} \sum_{n \in N \atop n < x^c} {x \over n}
\nonumber \\
&=&\ {1 \over x} \sum_{n \in N \atop n < x^c} \left (
\#\{m \leq x\ :\ n | m\} + O(1) \right ) \nonumber \\
&=&\ {1 \over x} \#\{ m \leq x\ :\ \exists n \in N,\ n < x^c,\ {\rm where\ } n|m\}
\ +\ O(x^{c-1}) \nonumber \\
&=&\ 1-\epsilon',
\end{eqnarray} 
where $0 < \epsilon' < 2\epsilon$ for $x$ sufficiently large.

The main result we will use in the development of our ``smooth sieve method''
is the following theorem:

\begin{theorem} \label{main_theorem}
For every $0 < \epsilon < \epsilon_0$ (for some $\epsilon_0$) 
and for $x$ sufficiently large
(in terms of $\epsilon$) there exists an LG set of integers 
$N \subseteq \{1,2,...,x\}$ with parameters $\epsilon$, $c=c(\epsilon)$
(that is, $c$ depends only on $\epsilon$, and not on $x$) and $x$.
Moreover, the following is an explicit example of
such a set:  For a certain constant $\delta$, depending only 
on $\epsilon$, and for $x$ sufficiently large, we let $N$ be the 
set of integers $n$ of the form
$p_1 p_2 \cdots p_k \leq x$ ($k$ variable), where the $p_i$'s are prime 
numbers such that 
$$
p_1 > p_2 > \cdots > p_k > x^\delta,
$$
and
\begin{equation} \label{prime_conditions}
{\rm For\ i=1,2,...,k-1},\ {x \over p_1p_2\cdots p_i} \geq p_i;\ \ 
{\rm and\ } 1 \leq {x \over p_1p_2\cdots p_k} < p_k.
\end{equation}
\end{theorem}

\noindent {\bf Remark 1:}  For the explicit construction given, in order
for condition 2 for being an LG set to be satisfied, we need only choose
$\delta$ so small that
$$
\sum_{n \in N} {1 \over n}\ >\ 1 - {\epsilon \over 2}.
$$ 
To see this, we note that since each member of $N$ has all its prime
divisors $> x^\delta$, basic sieve methods show that for $y > x^{1/2}$,
$$
\sum_{n \in N \atop n \in [y,2y]} {1 \over n}\ <\ {C_\delta \over \log x},
$$
where $C_\delta$ depends only on $\delta$.  Summing over dyadic intervals,
this implies
$$
\sum_{n \in N \atop x^c < n < x} {1 \over n}\ <\ D_\delta (1-c),
$$
for a certain constant $D_\delta$ depending only on $\delta$.  Now, by
taking $c$ sufficiently close to $1$, this sum can be made less than
$\epsilon/2$, which would give
$$
\sum_{n \in N \atop n < x^c} {1 \over n}\ >\ 1 - {\epsilon \over 2}
- \sum_{n \in N \atop x^c \leq n < x} {1 \over n}\ >\ 
1 - \epsilon.
$$ 
\bigskip

The main theorem of the paper is the following result:

\begin{theorem} \label{smooth_theorem}
Given $0 < \theta \leq 1$ and $\gamma > 0$, there exists 
$0 < \epsilon < 1$ so that for $x$ sufficiently large, if
$N$ is an LG set with parameters $\epsilon$, $c$, and $x$ as given in 
Theorem \ref{main_theorem}, then the following holds:
First, let $N_1$ be the set of $x^\theta$-smooths that lie in $N$
and are $< x^c$; let $N_2$ be the set of integers that are 
not $x^\theta$-smooth that lie in $N$ and are $<x^c$; 
let $w(n) \geq 0$ be some weighting function on positive integers 
$\leq x$; and let 
$$
\sigma\ =\ \sum_{s \leq x} w(s).
$$
Further, suppose that the following two inequalities hold
\begin{equation} \label{N1_inequality}
\sum_{q \in N_1} \sum_{s \leq x \atop q | s} w(s)\ >\ (1-\gamma)
\sigma \sum_{q \in N_1} {1 \over q};\ {\rm and}
\end{equation}

\begin{equation} \label{N2_inequality}
\sum_{q \in N_2} \sum_{s \leq x \atop q | s} w(s)\ >\ (1-\gamma) \sigma
\sum_{q \in N_2} {1 \over q}.
\end{equation}
Then, we have that
$$
\left | \sum_{s \leq x \atop s\ {\rm is\ }x^\theta-{\rm smooth}} 
w(s)\ -\ \sigma \sum_{q \in N_1} {1 \over q} \right |\ <\ 
2\gamma \sigma.
$$ 

Moreover, if one is only able to show (\ref{N1_inequality}), then one
can deduce the sharp lower bound
\begin{eqnarray}
\sum_{s \leq x \atop s\ {\rm is\ }x^\theta-{\rm smooth}} 
w(s)\ &>&\ (1-\gamma) \sigma \sum_{q \in N_1} {1 \over q} \nonumber \\
&>&\ (1-\gamma)(\rho(1/\theta) -2\epsilon) \sigma, \nonumber
\end{eqnarray}
for sufficiently large $x$.  
\end{theorem}

\noindent {\bf Remark.}  For a fixed $\epsilon > 0$, as $x$ tends to
infinity the sum of $1/q$ over $q \in N_1$ tends to 
a limit which is within $2\epsilon$ of $\rho(1/\theta)$.  So,
as $\epsilon$ tends to $0$ and $x$ tends to
infinity, the sum of $1/q$ 
over $q \in N_1$ tends to the limit $\rho(\theta^{-1})$, where $\rho$
is Dickman's function, and was stated earlier.  
\bigskip

Let us now discuss this theorem, to get a feel for how one might use it.
What this theorem is saying is that if one has a set of integers $S$,
and if one can produce only good lower bounds (upper bounds are not needed)
for how many elements in $S$ there are that are divisible by elements $q \in N_1$,
and by elements $q \in N_2$, then one can deduce that $S$ contains the expected number 
of $x^\theta$-smooths.   One of the strengths of the theorem is that the $q$'s for which 
one needs the divisibility conditions (\ref{N1_inequality}) and (\ref{N2_inequality}) to hold 
are less than a power of $x$ (where the power is less than $1$).  With a naive 
Buchstab identity approach one {\it must} confront the problem of whether there are
lots of elements in the set $S$ that are divisible by large primes; and, even if one
can solve that problem, there are still other problems involving divisibility by large
$q$'s that must be addressed.  These difficulties never need to be confronted in
the theorem above.

Another strength is that we only require lower bounds on the left-hand-side sums in
(\ref{N1_inequality}) and (\ref{N2_inequality}), and it is sometimes much easier to
produce such lower bounds, than it is to give asymptotic estimates.  
A fairly simple example is the following:  Suppose that one takes 
$\theta \in (0,1]$ and $\gamma > 0$ ``close'' to $0$, and lets $c=c(\theta,\gamma)$ be as 
in the theorem above.  Now take $A$ to be a subset of the integers $\leq x$
having at least $x^{c+\nu}$ elements, where $0 < \nu < 1$ can be taken arbitrarily
small.  Let $w(n)$ be the number of ways of writing $n$ as a difference of two elements
of $A$.  The sum of $w(n)$ over all positive integers $n$ is obviously 
${|A| \choose 2}$. 
It is also easy to see that the sum over all $w(n)$ with $n \geq 1$ 
divisible by $q$ is 
$$
\sum_{a=0}^{q-1} {A(a,q) \choose 2},
$$
where $A(a,q)$ is the number of elements of $A$ that are $\equiv a \pmod{q}$.
This expression is minimized if the elements of $A$ are as equidistributed amongst
the residue classes modulo $q$ as is possible; and so, this expression can be 
shown to be at least $\sim |A|^2 / (2q)$ in size.  
Summing over all $n$ divisible by
elements $q \in N_i$ of $w(n)$, we get that 
$$
\sum_{q \in N_i} \sum_{n \geq 1 \atop q | n} w(n)\ \gtrsim\ 
{|A|^2 \over 2} \sum_{q \in N_i} {1 \over q}.
$$
Thus, (\ref{N1_inequality}) and (\ref{N2_inequality}) hold for $x$ sufficiently large, 
and it follows that the number of
differences $a-b > 0$, $a,b \in A$, which are $x^\theta$-smooth is
``close'' to 
$U = (|A|^2/2) \sum_{q \in N_1} 1/q \approx \rho(1/\theta) |A|^2/2$; 
and, the smaller we take $\gamma$ to be, the closer this count will be to $U$.
In the next section, we will give a proof of a related (but more difficult) result.

The proof of theorem \ref{smooth_theorem} is so simple that we will give it here in
the introduction:
\bigskip

\noindent {\bf Proof of Theorem \ref{smooth_theorem}.}  
First, we want the value of $\epsilon$ for the LG set to be so small
that $\theta > \delta$, where $\delta$ is the parameter given in the
construction in Theorem \ref{main_theorem}.  There will be additional
demands on $\epsilon$ that we will give as the proof progresses.

Now suppose $n \leq x$ is divisible by some $q \in N$ (which is unique).  
Then,  we have that $n$ is 
$x^\theta$-smooth if and only if $q \in N_1$.  To see
this, first note that if $n = qk$, where $q \in N_1$,
then by the definition of the set $N$, we have that $x/q$ is less than
the smallest prime dividing $q$; that is, $k$ is less than $x^\theta$ 
(which must be greater than $x^\delta$).  
So, $n$ is $x^\theta$-smooth.  Conversely, suppose that
$n \leq x$ is $x^\theta$-smooth and $q | n$.  Then, if $q$ is not 
$x^\theta$-smooth, then neither is $n$.  

We deduce that the sum of $w(s)$ over all the elements $s \leq x$ 
that are $x^\theta$-smooths is {\it at least} the sum of $w(s)$ over all 
$s \leq x$ divisible by some $q \in N_1$.  This quantity is given
by the left hand side of (\ref{N1_inequality}).  On the other
hand, the sum of $w(s)$ over all $s \leq x$ that are $x^\theta$-smooth
is {\it at most} 
$\sigma - \tau$, where $\tau$ is the sum of $w(s)$ over all 
$s\leq x$ divisible by some $q \in N_2$.  From (\ref{N2_inequality}),
we get that this quantity is at most 
\begin{eqnarray}
\sigma\ -\ \tau\ &<&\ \sigma - (1-\gamma) \sigma \sum_{q \in N_2} {1 \over q}  
\nonumber \\
&<&\ \sigma - (1-\gamma) \sigma \left ( (1-\epsilon') -
\sum_{q \in N_1} {1 \over q} \right ) \nonumber \\
&<&\ \sigma \sum_{q \in N_1} {1 \over q}
\ +\ (\gamma + \epsilon') \sigma, \nonumber 
\end{eqnarray}
where $\epsilon'$ is as given in (\ref{N_sum}), and therefore depends on
the value of $\epsilon$ (but tends to $0$ as $\epsilon$ tends to $0$).
Now, if $\epsilon$ is sufficiently small, then $\epsilon'$ will be
smaller than $\gamma$, and so this last chain of inequalities would give
$$
\sigma\ -\ \tau\ <\ \sigma \sum_{q \in N_1} {1 \over q}\ +\ 2\gamma \sigma.
$$ 
Combining this with (\ref{N1_inequality}) then proves the theorem. 
$\blacksquare$
\bigskip

The remainder of this paper is organized as follows.  In the next section
we give an application of Theorem \ref{smooth_theorem} 
to counting the number of smooth sums
$a+b$, where $a$ lies in some set of integers $A$, and $b$ lies in
a set $B$, and in section \ref{main_theorem_section}, we will give 
a proof of Theorem \ref{main_theorem}. 
\bigskip

\section{An Application of Theorems \ref{smooth_theorem}}
\bigskip

Given sets of integers $A$ and $B$, which are subsets of 
$\{1,2,...,x\}$ having $\gg x$ elements each, 
it is an interesting and studied question to determine 
the number of $y$-smooth sums $a+b$, $a \in A$, $b \in B$.
There are several ways of attacking this sort of problem, one of which is to
use the circle method and exponential sums over smooth numbers, and
another is to use the large sieve.  We could also ask how 
$\tau(a+b)$ is distributed, or how large $P(a+b)$ can be,
where $\tau$ is the divisor function, and where $P(n)$ denotes the largest
prime factor of $n$.  Using the large sieve and the circle method,
these types of questions were given a thorough treatment in a series
of beautiful papers by A. Balog and A. Sarkozy \cite{balog1}, 
\cite{balog2}, \cite{balog3}, and \cite{balog4};  
P. Erd\H os, H. Maier, and A. S\'ark\H ozy \cite{erdos}; 
A. S\'ark\H ozy and C. L. Stewart \cite{sarkozy1},
\cite{sarkozy2}, \cite{sarkozy3}, \cite{sarkozy4}; 
C. Pomerance, A. S\'ark\H ozy, and C. L. Stewart \cite{pomerance};
and R. de la Bret\`eche \cite{breteche}.
The paper by de la Br\`eteche is more relevant to
the main result of this section, and we give here one of his theorems:

\begin{theorem}
Suppose that $A$ and $B$ are subsets of the integers in $\{1,2,...,x\}$.
For a given integer $y \leq x$, let $u = (\log x)/\log y$.  Then,
uniformly for $x \geq 3$, $\exp ((\log x)^{2/3+\epsilon}) < y \leq x$
we have
\begin{eqnarray}
&& \#\{a \in A,\ b \in B\ :\ P(a+b) \leq y\}\nonumber \\  
&&\hskip1in =\ 
|A|\cdot |B| \rho(u) \left ( 1 + O \left ( {x \over \sqrt{|A| \cdot |B|}}
 {\log (u+1) \over \log y} \right ) \right ), \nonumber
\end{eqnarray}
where $P(n)$ denotes the largest prime divisor of $n$, and 
where $\rho$ is Dickman's function.
\end{theorem} 

R. de la Bret\`eche used estimates for exponential sums and the circle
method to prove this result.  Notice that if 
$|A|\cdot |B| \ll (x/\log x)^2$, then his result fails to
prove that there are the expected number of sums that are $y$-smooth
for any $y < x$, because in this case the big-Oh term is
$\gg 1$.  

Let us now consider what happens in the case when $y = x^\theta$ 
(and so $u = 1/\theta$):  Is it possible to show that if 
$|A|, |B| > x^c$, for some $0 < c < 1$, then we get the expected
number of sums $a+b$ being $y$-smooth?  It is easy to see that the
answer is no, no matter how close to $1$ we take $c$ to be.
For example, we could take $A$ and $B$ to both be
the set of integers $\leq x$ that are divisible by some prime
number $p$ around size $x^{1-c}$.  Notice here that $|A| \sim x^c$.
But then, the sums $a+b$, $a,b \in A$ are numbers of the form 
$p k$, where $k < 2x^c$, and such a sum is $x^\theta$-smooth 
if and only if $k$ is $x^\theta$-smooth, when $p < x^\theta$.  
Thus, one would expect (and can show) that
$$
{\#\{a,b \in A\ :\ a+b\ {\rm is\ }x^\theta-{\rm smooth}\} \over
|A|\cdot |B|}\ \sim\ {\Psi(2x^c, x^\theta) \over 2x^c}\ \sim\ 
\rho(c/\theta).
$$ 
On the other hand, 
the proportion of $x^\theta$ smooths $\leq x$ is $\sim \rho(u')$,
where $u' = \log(x)/\log(x^\theta) = 1/\theta$, which is not
$c/\theta$.  So, the type of result we might try to prove is the
following:

\begin{theorem} \label{smooth_application}  
Given $0 < \theta \leq 1$, and
$\gamma > 0$, there exists $\nu = \nu(\theta,\gamma) \in (0,1)$ 
so that the following
holds:  For $x$ sufficiently large, if 
$A,B \subseteq \{1,2,...,x\}$ satisfy $|A|, |B| > x^\nu$, then
\begin{eqnarray}
&& \left | \#\{a, b\ :\ a\in A,\ b \in B,\ 
a+b\ {\rm\ is\ }x^\theta-{\rm smooth}\}
- \rho(1/\theta) |A|\cdot |B| \right |\nonumber \\
&&\hskip1in <\ \gamma |A|\cdot |B|. \nonumber
\end{eqnarray}
The same result holds for differences $a-b$ (with possibly a different
value for $\nu$).
\end{theorem}

There are other methods for proving this type of theorem, besides the
exponential sums approach used by de la Bret\`eche, 
such as an application of Buchstab's
identity, together with a form of the large sieve for sieving by
composite moduli; however, these methods are technical, and it does
not seem possible to give an easy and elegant proof of Theorem \ref{smooth_application}
using them.  

In the remainder of this section we will give an entirely elementary
proof of Theorem \ref{smooth_application} using Theorem \ref{smooth_theorem}.
First, though we need a ``large sieve''-type result for composite moduli,
and although such a result can be easily proved by simply modifying the 
standard large sieve, we give here (perhaps astonishingly) a completely 
elementary proof based on LG sets.

\begin{theorem} \label{local_global_sieve}
Given $\epsilon > 0$ and $x$ sufficiently large, let $N$ be 
an LG set for parameters $\epsilon$ and $x$, as given in 
Theorem \ref{main_theorem}.  Further, suppose that $c=c(\epsilon)$
is as in property 2 for being an LG set.  
Suppose that $C \subseteq \{1,2,...,x\}$, and let $C(a,q)$
denote the number of elements of $C$ that are congruent to 
$a$ modulo $q$.  Then, we have that
$$
\sum_{q \in N \atop q < x^c} \sum_{a=0}^{q-1} \left ( 
C(a,q) -\ {|C| \over q} \right )^2\ <\ |C|(2\epsilon |C| + x^c).
$$  
\end{theorem}

The proof of this ``Large Sieve''-like theorem is so simple, we will
not postpone its proof to a later section:
\bigskip

\noindent {\bf Proof.  }  
To make the notation simple, when we sum over $q \in N$, we mean
the sum over those $q \in N$ satisfying $q \leq x^c$.

We note that if $b,c \in C$, and $b \neq c$, then if 
$q \in N$ divides $b-c$, we must have that $q$ is unique;
otherwise, if $q' \in N$ also divides $b-c$, then 
lcm$(q,q') > x$ divides $b-c$, which is impossible.

Thus, we have that
\begin{eqnarray}
|C|^2\ &>&\ 
\sum_{q \in N} \#\{b,c \in C,\ b \neq c\ :\ q | (b-c)\} \nonumber \\
&=&\ \sum_{q \in N} \sum_{a=0}^{q-1} \left ( C(a,q)^2 - C(a,q) \right ) 
\nonumber \\
&=&\ \sum_{q \in N} \sum_{a=0}^{q-1} C(a,q)^2\ -\ x^c|C|. \nonumber
\end{eqnarray}
Thus,
$$
|C| (x^c + |C|)\ >\ \sum_{q \in N} \sum_{a=0}^{q-1} C(a,q)^2.
$$

It now follows that
\begin{eqnarray}
\sum_{q \in N} \sum_{a=0}^{q-1} 
\left ( C(a,q) - {|C|\over q} \right )^2
\ &=&\ \sum_{q \in N} \sum_{a=0}^{q-1} C(a,q)^2 - |C|^2 \sum_{q \in N} 
{1 \over q}\nonumber \\
&<&\ (1 - (1-\epsilon'))|C|^2 + x^c |C| \nonumber \\
&<&\ |C| (2\epsilon |C| + x^c), \nonumber
\end{eqnarray}
where $\epsilon'$ is as in (\ref{N_sum}).  The theorem is now
proved.\ \ \ \ \ \ \ $\blacksquare$
\bigskip

\noindent {\bf Proof of Theorem \ref{smooth_application}.}  
Given $0 < \theta \leq 1$, and $\gamma > 0$, we suppose $\epsilon$
is so small that the conclusion of
Theorem \ref{smooth_theorem} holds for $x$ sufficiently large.
Let $N$ be the LG set with parameters 
$\epsilon$, $c=c(\epsilon)$ and $x$ (for $x$ sufficiently large) as appears in the construction
in Theorem \ref{main_theorem}.  Finally, let $\delta = \delta(\epsilon)$ 
be the parameter also given in this construction,
and let $N_1$ and $N_2$ be the sets as described in the 
conclusion of Theorem \ref{smooth_theorem}.  Two additional demands on
$\epsilon$ and $x$ is that we will need $\epsilon$ to be so small,
and $x$ so large that
\begin{equation} \label{theta_gamma_condition}
\left | \sum_{n \in N_1} {1 \over n}\ -\ \rho(1/\theta) \right |
\ <\ {\gamma \over 4},
\end{equation}
and
$$
\epsilon\ <\ {\gamma \over 12} \sum_{q \in M} {1 \over q},\ {\rm for\ }
 M = N_1\ {\rm and\ }N_2.
$$ 

We will show that the conclusion of our theorem holds for 
any $\nu > c = c(\epsilon)$ for $x$ sufficiently large.

Let $\alpha$ be the indicator function on the set $A$,
let $\beta$ be the indicator function on the set $B$,
let 
$$
A(a,q)\ =\ \sum_{n\leq x \atop n \equiv a \pmod{q}} \alpha(n),\ {\rm and\ }
B(a,q)\ =\ \sum_{n \leq x \atop n \equiv a \pmod{q}} \beta(n),
$$
and, finally, define the weight function
$$
w(n)\ =\ \sum_{a + b = n} \alpha(a) \beta(b).
$$
Then, for $M=N_1$ or $N_2$, we get 
\begin{eqnarray} \label{w_focus}
\sum_{q \in M} \sum_{n \leq x \atop q | n} w(n)\ &=&\ 
\sum_{q \in M} \sum_{a=0}^{q-1} A(a,q) B(q-a,q) \nonumber \\
&=&\ \sum_{q \in M} \sum_{a=0}^{q-1} \left (A(a,q) - {|A| \over q}
\right ) \left ( B(q-a,q) - {|B| \over q} \right ) \nonumber \\
&&\ \ \ \ +\ |A|\cdot |B| \sum_{q \in M} {1 \over q}. 
\end{eqnarray}

Now, then, to bound the last double sum in (\ref{w_focus}) 
from above we apply the Cauchy-Schwarz inequality, together with
Theorem \ref{local_global_sieve} with $C=A$ and $C=B$:
\begin{eqnarray} \label{error_cauchy}
&& \sum_{q \in M} \sum_{a=0}^{q-1} \left | A(a,q) - {|A| \over q}
\right |\ \left | B(q-a,q) - {|B| \over q} \right |
\nonumber \\
&&\ \ \ \leq\ \left ( \sum_{q \in M} \sum_{a=0}^{q-1} 
\left ( A(a,q) - {|A| \over q} \right )^2 \right )^{1/2}
\left ( \sum_{q \in M} \sum_{a=0}^{q-1} \left ( B(a,q) - {|B| \over q}
\right )^2 \right )^{1/2} \nonumber \\
&&\ \ \ \leq\ (|A|(2 \epsilon |A| + x^c) )^{1/2} (|B|(2\epsilon |B| + x^c)^{1/2}
\nonumber \\
&&\ \ \ \leq\ 3\epsilon |A|\cdot |B|,
\end{eqnarray}
for $|A|, |B| > x^c/\epsilon$, which certainly holds if $\nu > c$ and $x$
is sufficiently large.

Combining this with (\ref{w_focus}) we deduce that
\begin{eqnarray}
\sum_{q \in M} \sum_{n \leq x \atop q|n} w(n)\ &>&\ |A|\cdot |B|
\left ( -3\epsilon +  \sum_{q \in M} {1 \over q}\right ) \nonumber \\
&>&\ |A|\cdot |B| \left ( 1 - {\gamma \over 4} \right ) \sum_{q \in M}
{1 \over q}.
\end{eqnarray}
for $M = N_1$ or $N_2$.  Thus, the conditions of Theorem \ref{smooth_theorem}
are met, and we deduce that if $\sigma$ is the sum of $w(n)$ over all
$n \leq x$, which is $|A|\cdot |B|$, then
\begin{eqnarray}
&& | \#\{ a,b\ :\ a \in A,\ b \in B,\ a+b\ {\rm is\ }x^\theta-{\rm smooth}\}
- \rho(1/\theta) \sigma | \nonumber \\ 
&& \leq\ 
\left | \#\{ a,b\ :\ a \in A,\ b \in B,\ a+b\ {\rm is\ }x^\theta-{\rm smooth}\}
- \sigma \sum_{q \in N_1} {1 \over q} \right | + {\gamma \over 4} \sigma 
\nonumber \\
&&<\ {\gamma \over 2} \sigma + {\gamma \over 4} \sigma
\nonumber \\
&&<\ \gamma \sigma. \nonumber
\end{eqnarray}
The theorem now follows.  $\blacksquare$

\section{Proof of Theorem \ref{main_theorem}} \label{main_theorem_section}
\bigskip

First, we show that the set $N$ described in the statement of the theorem
satisfies the first condition for being an LG set, namely that for
any distinct pair of integers $n_1, n_2 \in N$, we have lcm$(n_1,n_2) > x$:
Given such $n_1,n_2$, write out their prime factorizations as
\begin{eqnarray}
n_1\ &=&\ p_1 \cdots p_k,\ p_1 > p_2 > \cdots > p_k;\ {\rm and} \nonumber \\
n_2\ &=&\ q_1 \cdots q_\ell,\ q_1 > q_2 > \cdots > q_\ell. \nonumber
\end{eqnarray}
Without loss of generality, we can assume that $p_k \leq q_\ell$.

Now, if there is some prime $q_i$ which is distinct from the primes
$p_1,...,p_k$, then we would have that the lcm of $n_1$ and $n_2$
is divisible by the product of primes $q_i p_1\cdots p_k$, and this product
exceeds $x$, because from (\ref{prime_conditions})
$$
q_i p_1 \cdots p_k\ >\ q_i {x \over p_k}\ \geq\ q_\ell {x \over p_k}\ \geq\ x.
$$

So we are left to consider what happens when the $q_i$'s are a subset of the 
$p_i$'s.  We break this case into two sub-cases, with the first one where
$p_k = q_\ell$, and the second where $p_k < q_\ell$.  

In the case 
$p_k = q_\ell$, we must have that there exists one of the primes 
$p_i > q_\ell$ such that $p_i$ is distinct from $q_1,...,q_\ell$,
since otherwise we would have $n_1 = n_2$.  But now our assumption 
gives
$$
{x \over p_1 \cdots p_k}\ \leq\ {x \over p_i q_1 \cdots q_\ell}
\ <\ {q_\ell \over p_i}\ <\ 1,
$$
which is impossible.

So, we may assume $p_k < q_\ell$.  For this case, let $j < k$ be the 
index where $p_j = q_\ell$ (which exists since $q_i$'s are a subset of 
the $p_i$'s).  Then, we have 
$$
p_1\cdots p_j\ \geq\ q_1\cdots q_\ell.
$$
From (\ref{prime_conditions}) this gives
$$
q_\ell\ =\ p_j\ \leq\ 
{x \over p_1\cdots p_j}\ \leq\ {x \over q_1 \cdots q_\ell}\ 
<\ q_\ell,
$$
which is impossible.  So, we conclude that the set described in the
statement of the theorem satisfies the first condition for being
an LG set.
\bigskip

We are left to show that for some $0 < \delta < 1$, and all 
$x$ sufficiently large, all but $\epsilon x$ of the integers $m \leq x$
are divisible by some member of $N$.  We will do this by identifying
a subset $T \subseteq \{1,2,...,x\}$ having at least $(1-\epsilon/2)x$
elements, such that all but at most $\epsilon x/2$ of the elements of $T$
are divisible by some member of $N$.  The way we will show this is to
construct a weighting function $f(t) > 0$ such that if 
$t \in T$ is not divisible by any $n \in N$, then $f(t)$ will be
``large''.  But then, we will show that average value of $f(t)$
over all $t \in T$ is ``much smaller'' these large values; so,
it will follow that there can be few integers $t \in T$ not 
divisible by any $n \in N$.
\bigskip

Let $k = \lfloor 1/\epsilon \rfloor + 1$; let $0 < \gamma < 1$ be 
some constant to be chosen later; let $H_1,...,H_k$ be intervals
given by $H_j = [x^{\gamma^{j+1}}, x^{\gamma^j})$; let $I_1,...,I_k$
be the intervals where $I_j = [x^{\gamma^{j+1}}, x^{\gamma^j/2}]$;
and finally, let $J_1,...,J_k$ be the intervals 
$J_j = (x^{\gamma^j/2},x^{\gamma^j})$.  Note here that
$H_j = I_j \cup J_j$.  Our constant $\delta$ in
the construction in the statement of the main theorem will 
be $\delta = \gamma^{k+1}$.  

We first claim that for $x$ sufficiently large, all but 
at most $\epsilon x/4$ integers $m \leq x$ satisfy the following 
inequality for all $j=1,2,...,k$:
\begin{equation} \label{small_sum_equation}
\sum_{p^a|m,\ p\ {\rm prime} \atop p < x^{\gamma^j}} \log p\ <\ 
5k^2 \gamma^j \log x,
\end{equation}
To see this, we first note that for any $j=1,2,...,k$, 
\begin{eqnarray}\label{gamma_sum}
\sum_{m \leq x} \sum_{p^a | m,\ p\ {\rm prime} \atop p < x^{\gamma^j}} \log p
\ &=&\ \sum_{p \leq x^{\gamma^j} \atop p\ {\rm prime}} (\log p)
\sum_{a \geq 1} \#\{m \leq x\ :\ p^a | m\} \nonumber \\
&\leq&\ \sum_{p \leq x^{\gamma^j} \atop p\ {\rm prime}} (\log p) 
\sum_{a \geq 1} {x \over p^a} \nonumber \\
&=&\ x\sum_{p \leq x^{\gamma^j} \atop p\ {\rm prime}} 
{\log p \over p}\ +\ O (x) \nonumber \\
&=&\ \gamma^j x \log x\ +\ O(x), 
\end{eqnarray}
where the constant in the last big-oh depends on $\gamma$.  The last line
here was gotten by using the fact that the sum over primes $p \leq x$ of
$(\log p)/p$ is $\log x + O(1)$.  Now, for $x$ sufficiently large, 
there can be at most $x/(4k^2)$ integers $m \leq x$ 
which fail to satisfy (\ref{small_sum_equation}) for $j$, since otherwise
these exceptional integers $n$ would force the first sum in 
(\ref{gamma_sum}) to be of size at least $5\gamma^j x (\log x)/4$, and we
know this sum is of size at most $\gamma^j x (\log x)$.  So, for $x$
sufficiently large, (\ref{small_sum_equation}) holds for all
$j=1,2,...,k$ for all but at most
$$
k {x \over 4k^2}\ =\ {x \over 4k}\ <\ {\epsilon x \over 4}
$$  
exceptional integers $m \leq x$.  Let $S$ denote the set of integers satisfying
(\ref{small_sum_equation}) for all $j=1,2,...,k$.  Then, we have 
shown $|S| > (1-\epsilon/4)x$ for $x$ sufficiently large.

Let $h(s)$ be the number of integers $j=1,2,...,k$
such that $s$ is divisible by some prime $p \in J_j$.  
So, $0 \leq h(s) \leq k$.  We will show below that all but at most 
$\epsilon x/4$ integers $s \leq x$ satisfy 
\begin{equation} \label{T_inequalities}
\left | h(s) - {k \over 2} \right | <  k^{2/3};\  
s > {\epsilon x \over 100};\ 
{\rm and,\ } p^2 | s,\ p\ {\rm prime}\ \Rightarrow\ p < x^\delta.
\end{equation}
These last two conditions are obviously satisfied for all but at most 
$\epsilon x/50$ integers $s \leq x$ for large $x$; and so, we just 
need to show that the first condition holds for all but at most $\epsilon x/5$ 
integers $n \leq x$.  We will then let $T$ be the set of all $s \in S$ satisfying these
conditions.  Clearly we will have $|T| > (1-\epsilon/2)x$.
\bigskip

To prove that the first condition of (\ref{T_inequalities}) holds for all but at most
$\epsilon x/5$ integers $s \leq x$, we begin by supposing $V \subseteq \{1,2,...,k\}$ and
$W = \{1,2,...,k\} \setminus V$, and letting
$U$ denote the set of all integers $u \leq x$ such that
\bigskip

1.  For every $v \in V$, $u$ is not divisible by any prime $p \in J_v$; and,

2.  For every $w \in W$, $u$ is divisible by some prime $p \in J_w$.
\bigskip
 
\noindent To estimate the size of $|U|$ we require the following corollary
of the combinatorial sieve of Rosser:

\begin{proposition}
For every $0 < \beta < 1$, there exists $0 < \tau < 1$ 
so that for $x$ sufficiently large, the following holds:
Suppose $P$ is a subset of the primes $\leq x^\tau$, and let 
$Z$ be the set of all integers $\leq x$ not divisible by 
any prime $p \in P$.  Then, if we let
$$
\Delta\ =\ \prod_{p \in P} \left ( 1 - {1 \over p} \right ),
$$
we will have
$$
x(1-\beta)\Delta\ <\ |Z|\ <\ x(1+\beta)\Delta.
$$
\end{proposition}

Now, to estimate the size of $|U|$, for each indexing set 
$V' \supseteq V$, we let $Z(V')$ denote the set of all integers
$z \leq x$ not divisible by any prime $p \in J_v$, for any 
$v \in V'$.  Then, by a simple inclusion-exclusion argument
we have
$$
|U|\ =\ \sum_{V' \supseteq V} (-1)^{|V'|-|V|} |Z(V')|.
$$ 
In order to apply the proposition, we first estimate
\begin{eqnarray}
\Delta\ &=&\ \prod_{v \in V'} \prod_{p \in J_v \atop p\ {\rm prime}} 
\left ( 1 - {1 \over p} \right ) \nonumber \\
&=&\ \prod_{v \in V'} \exp \left ( - \sum_{p \in J_v \atop p\ {\rm prime}}
{1 \over p}\ +\ O \left ( {1 \over x^{\gamma^j/2}} \right ) \right ) \nonumber \\
&=&\ \prod_{v \in V'} \exp \left ( -\log 2 + O \left ( {1 \over \gamma^j\log x}
\right ) \right ) \nonumber \\
&=&\ {1 \over 2^{|V'| + O(k/\gamma^j \log x)}}. \nonumber 
\end{eqnarray}
Note that this last expression is asymptotically $2^{-|V'|}$.

It now follows that from the proposition above, and from our equation
for $|U|$ that for fixed $\epsilon, \beta$, if $\gamma$ is sufficiently
small and $x$ is sufficiently large, then
\begin{eqnarray}
|U|\ &=&\ x \left ( 1 + O \left ( {1 \over \log x} \right ) \right )
\sum_{V' \supseteq V} (-1)^{|V'| - |V|} {1 \over 2^{|V'|}}
\ +\ O \left ( 2^k \beta x \right ) \nonumber \\
&=&\ {x \over 2^{|V|}} \left ( 1 + O \left ( {1 \over \log x} \right ) \right )
\sum_{V'' \subseteq \{1,2,...,k\} \setminus V}
{(-1)^{|V''|} \over 2^{|V''|}}\ +\ O \left ( 2^k \beta x \right ) \nonumber \\
&=&\ {x \over 2^k}\ +\ O \left ( 2^k \beta x \right ). \nonumber 
\end{eqnarray}
The implied constant in this last big-oh depends on $\gamma$ and $\epsilon$.

So, the number of integers $s \leq x$ such that $h(s)$ is more than
$k^{2/3}$ away from $k/2$ is
\begin{eqnarray}
\sum_{V \subseteq \{1,2,...,k\} \atop ||V| - k/2|\ >\ k^{2/3}} 
\left ( {x \over 2^k} + O \left ( 2^k \beta x \right ) \right ) 
\ &=&\ {x \over 2^k} \sum_{0 \leq j \leq k \atop 
|j-k/2| > k^{2/3}} {k \choose j}
\ +\ O \left ( 4^k \beta x \right ). \nonumber
\end{eqnarray}
Now, the sum over these binomial coefficients can be shown to be 
smaller than $2^k/k^2$; and so, for $\epsilon_0 < 1/5$,
by choosing $\beta$ and $\gamma$ sufficiently small, we will have 
for $x$ sufficiently large that there can be at most
$$
{x \over k^2} + O\left ( 4^k \beta x \right )\ <\ 
{\epsilon x \over 5}
$$
integers $s \leq x$ with $|h(s) - k/2| > k^{2/3}$.  This then proves that 
(\ref{T_inequalities}) holds for all but at most $\epsilon x/4$ elements
of $S$.
\bigskip

To finish the proof of the theorem we will focus on the following
function:
$$
f(s)\ =\ \prod_{j=1}^k \left ( \gamma^{2j} \log^2 x
\ +\ \left ( \sum_{p \in I_j,\ p|s \atop p\ {\rm prime}} \log p 
\right )^2 \right ).
$$
First, let us calculate the sum of $f(s)$ over all $s \leq x$:
This sum equals
$$
\sum_{V \subseteq \{1,2,...,k\}} (\log x)^{2|V|} 
\prod_{v \in V} \gamma^{2v} \sum_{s \leq x} 
\prod_{w \in \{1,2,...,k\} \setminus V}
\left ( \sum_{p \in I_w,\ p|s \atop p\ {\rm prime}} \log p
\right )^2. 
$$
We now work with this inner sum by fixing a $V \subseteq \{1,2,...,k\}$,
and letting $W = \{1,2,...,k\} \setminus V$, and then this inner sum 
is
\begin{eqnarray}
&&\sum_{s \leq x} \prod_{w \in W} \left ( \sum_{p \in I_w,\ p|s
\atop p\ {\rm prime}} \log p \right )^2 \nonumber \\
&&\ \ \ \ \ \leq\ x \prod_{w \in W} 
\left ( 2\sum_{p_1 < p_2;\ p_1,p_2 \in I_w \atop p_1, p_2\ {\rm prime}}
{(\log p_1)(\log p_2) \over p_1p_2}\ +\ 
\sum_{p \in I_w \atop p\ {\rm prime}} {\log^2 p \over p} \right ) \nonumber \\
&&\ \ \ \ \ =\ x \prod_{w \in W} 
\left ( \left ( \sum_{p \in I_w \atop p\ {\rm prime}} {\log p \over p}
\right )^2 + \sum_{p \in I_w \atop p\ {\rm prime}}
{\log^2 p \over p} \right ) \nonumber \\
&&\ \ \ \ \ \leq\ x \prod_{w \in W}
\left ( \left ( \sum_{p \leq x^{\gamma^w/2} \atop p\ {\rm prime}} {\log p \over p}
\right )^2 + \sum_{p \leq x^{\gamma^w/2} \atop p\ {\rm prime}} 
{\log^2 p \over p} \right ) \nonumber \\
&&\ \ \ \ \ =\ x \prod_{w \in W} 
\left ( {\gamma^{2w} \over 4} \log^2 x + {\gamma^{2w} \over 8} \log^2 x 
+ O ( \gamma^w \log x ) \right ) \nonumber \\
&&\ \ \ \ \ =\ x \left ( 1 + O \left ( {1 \over \log x} \right ) \right )
\left ( {3 \log^2 x \over 8} \right )^{|W|} 
\gamma^{2 \sum_{w \in W} w}, \nonumber
\end{eqnarray} 
where the implied constant in this last big-oh depends on $\gamma$ and
$|W| \leq k$.  

We then get that the sum of $f(s)$ over $s \leq x$ is at most
\begin{eqnarray}
&& \gamma^{k(k+1)} x \left ( 1 + O \left ( {1 \over \log x} \right ) \right )
(\log x)^{2k} \sum_{V \subseteq \{1,2,...,k\}}
\left ( {3 \over 8} \right )^{k-|V|} \nonumber \\
&& =\ \gamma^{k(k+1)} x \left ( 1 + O \left ( {1 \over \log x} \right ) \right )
(\log x)^{2k} \left ( {11 \over 8} \right )^k. \nonumber 
\end{eqnarray}
\bigskip

Before we bound $f(t)$ from below for an arbitrary $t \in T$ that
fails to be divisible by any $n \in N$, we make a general observation:
If $p \geq x^\delta$ is any prime divisor of such a $t$,
then we must have that 
$$
{x \over \prod_{q \geq p,\ q|t \atop q\ {\rm prime}} q}\ >\ p,
$$
for otherwise $t$ is divisible by some integer $n \in N$. 
Now, since each member of $t \in T$ is at least $\epsilon x/100$, 
we deduce that
\begin{eqnarray} \label{product_inequality}
\prod_{q \leq p,\ q^a||t \atop q\ {\rm prime}} q^a
\ &=&\ {t \over \prod_{q \geq p,\ q | t \atop q\ {\rm prime}} q}
\nonumber \\
&\geq&\ {\epsilon x \over 100 \prod_{q \geq p,\ q | t \atop 
q\ {\rm prime}} q} \nonumber \\
&>&\ {\epsilon p \over 100}. 
\end{eqnarray} 
Now let $j$ be one of the $> k/2 - k^{2/3}$ indices for which 
$t$ is not divisible by any prime $p \in J_j$.  Further, suppose that
$j$ is not the smallest index, which guarantees that $t$ is divisible 
by at least some prime greater than $x^{\gamma^j}$.  
Then, the fact that (\ref{product_inequality}) must hold implies
$$
\prod_{q \leq x^{\gamma^j/2},\ q^a || t \atop q\ {\rm prime}} q^a
\ =\ \prod_{q \leq x^{\gamma^j},\ q^a || t \atop q\ {\rm prime}} q^a
\ >\ {\epsilon x^{\gamma^j} \over 100}.
$$
Taking logs of both sides gives
$$
\sum_{q \leq x^{\gamma^j/2},\ q^a | t \atop q\ {\rm prime}} 
\log q\ >\ \gamma^j \log x\ + \log(\epsilon/100) 
$$
Now, from our assumption (\ref{small_sum_equation}) we then deduce
\begin{eqnarray}
\sum_{q \in I_j,\ q | t \atop q\ {\rm prime}} \log q
\ &=&\ \sum_{q \leq x^{\gamma^j/2},\ q^a |t \atop q\ {\rm prime}}
\log q\ -\ \sum_{q \leq x^{\gamma^{j+1}},\ q^a | t \atop q\ {\rm prime}}
\log q \nonumber \\
&>&\ \gamma^j (1 - 5\gamma k^2) \log x\ +\ \log(\epsilon/100) \nonumber \\
&>&\ \gamma^j (1-\epsilon/2) \log x,
\end{eqnarray}
for a fixed $\epsilon$ and $\gamma$ sufficiently small.
Thus, if we let $X$ denote the set of indices $j$ such that
$t$ is not divisible by any prime $p \in J_j$, then 
\begin{eqnarray}
f(t)\ &\geq&\ \prod_{j \in X} \left ( \gamma^{2j}\log^2 x + 
(1-\epsilon/2)^2 \gamma^{2j} \log^2 x \right ) \prod_{j \in \{1,2,...,k\}
\setminus X} \gamma^{2j} \log^2x \nonumber \\
&\geq&\ \gamma^{k(k+1)} (\log^{2k} x) (2 - \epsilon)^{k/2 - k^{2/3}}.
\nonumber
\end{eqnarray}

So, if we let $Y$ be the integers $t \in T$ not divisible by any
$n \in N$, then we have
\begin{eqnarray}
|Y| \gamma^{k(k+1)} (\log^{2k} x) (2-\epsilon)^{k/2-k^{2/3}}
\ &<&\ \sum_{y\in Y} f(y)\ <\ \sum_{s \leq x} f(s) \nonumber \\
&<&\ x \gamma^{k(k+1)} (\log^{2k} x) 
\left ( {11 \over 8} \right )^k\nonumber \\
&&\hskip1in +\ O (x \log^{2k-1} x). \nonumber
\end{eqnarray}

So, if $\epsilon_0$ were sufficiently small (so as to make $k$ sufficiently
large, and $\epsilon$ sufficiently small), then we would clearly have
$|Y| < \epsilon x / 2$.  Therefore, all but $\epsilon x/2$ integers $t \in T$
is divisible by some $n \in N$ once $\epsilon_0$ is sufficiently small; 
and therefore, all but at most $\epsilon x$ integers $\leq x$ are divisible
by some $n \in N$.\ \ \ \ \ \ \ \ \ $\blacksquare$


\begin{thebibliography}{999}

\bibitem{baker} R. C. Baker and G. Harman, {\it Shifted Primes without Large Prime
Factors} Acta Arith. {\bf 83} (1998), 331-361.

\bibitem{balog1} A. Balog and A. S\'ark\H ozy, {\em On sums of Integers
Having Small Prime Factors. I, II}, Studia Sci. Math. Hungar. {\bf 19}
(1984), 35-47.

\bibitem{balog2} ----------------------------, {\em On Sums of Sequences
of Integers. I}, Acta Arith. {\bf 44} (1984), 73-86.

\bibitem{balog3} ----------------------------, {\em On Sums of Sequences
of Integers. II}, Acta Math. Hungar. {\bf 44} (1984), 169-179.

\bibitem{balog4} ----------------------------, {\em On Sums of Sequences
of Integers. III}, Acta Math. Hungar. {\bf 44} (1984), 339-349.

\bibitem{breteche} R. de la Bret\`eche, {\em Sommes sans Grand Facteur
Premier [Sums without Large Prime Factors]} Acta Arith. {\bf 88} (1999), 
1-14.

\bibitem{dickman} Dickman, {\it On the Frequency of Numbers Containing Prime
Factors of a Certain Relative Magnitude} Ark. Mat. Astr. Fys {\bf 22} (1930), 1-14.

\bibitem{erdos} P. Erd\H os, H. Maier, and A. S\'ark\H ozy,
{\em On the Distribution of the Number of Prime Factors of Sums
$a+b$} Trans. Amer. Math. Soc. {\bf 302} (1987), 269-280.

\bibitem{friedlander} J. B. Friedlander, {\it Shifted Primes without Large Prime
Factors}, NATO Adv. Sci. Inst. Ser. C Math. Phys. Sci., 265, Kluwer Acad. Publ.,
Dordrecht, 1989.

\bibitem{pomerance} C. Pomerance, A. S\'ark\H ozy, and C. L. Stewart,
{\em On Divisors of Sums of Integers. III.} Pacific J. Math. 
{\bf 133} (1988), 363-379.

\bibitem{sarkozy1} A. S\'ark\H ozy and C. L. Stewart, {\em On Divisors
of Sums of Integers. I} Acta Math. Hungar. {\bf 48} (1986), 147-154.

\bibitem{sarkozy2} ---------------------------------, {\em On Divisors
of Sums of Integers. II.}, J. Reine Angew. Math. {\bf 365} (1986), 171-191.

\bibitem{sarkozy3} A. S\'ark\H ozy and C. L. Stewart, {\em On Divisors
of Sums of Integers. IV.} Canad. J. Math. {\bf 40} (1988), 788-816.

\bibitem{sarkozy4} ---------------------------------, {\em On Divisors
of Sums of Integers. V.} Pacific J. Math. {\bf 166} (1994), 373-384.
 
\end{thebibliography}
\end{document}